\documentclass[12pt]{article}

\usepackage{typearea}
\usepackage{a4wide}
\usepackage{amsmath}
\usepackage{amsfonts}
\usepackage[frenchb, english]{babel}
\usepackage[T1]{fontenc} 
\usepackage[dvips,pdftex
]{graphicx}
\usepackage{mathrsfs}
\usepackage{cases}
\usepackage{aeguill}
\usepackage{eufrak}
\def\HH{\EuFrak H}

\newtheorem{theo}{Theorem}[section]
\newtheorem{prop}{Proposition}[section]
\newtheorem{corol}{Corollary}[section]

\newtheorem{lemme}{Lemma}[section]
\newtheorem{Rem}{Remark}[section]

\newcommand{\CQFD}{\hfill $\square$\\}

\newcommand{\ind}{\mathbf{1}}

\def\rit{\Bbb{R}}
\def\cit{\Bbb{C}}
\def\nit{\Bbb{N}}

\def\ee{\Bbb{E}}

\def\cvfree{\stackrel{\mbox{\tiny free}}{\longrightarrow}}
\def\cvlaw{\stackrel{\mbox{\tiny law}}{\longrightarrow}}
\def\loifree{\stackrel{\mbox{\tiny free}}{\sim}}

\def\pp{\Bbb{P}}
\def\E{\mathop{\hbox{\rm I\kern-0.20em E}}\nolimits}
\def\Var{\mathop{\hbox{\rm Var}}\nolimits}

\def\lto{\longrightarrow}

\def\og{\leavevmode\raise.3ex
     \hbox{$\scriptscriptstyle\langle\!\langle$~}}
\def\fg{\leavevmode\raise.3ex
     \hbox{~$\!\scriptscriptstyle\,\rangle\!\rangle$}~}

\allowdisplaybreaks

\begin{document}

\renewcommand{\thefootnote}{\fnsymbol{footnote}}

\renewcommand{\thefootnote}{\fnsymbol{footnote}}
\begin{center}
{\large{ \textsc{Asymptotic Cram\'er type decomposition for Wiener and Wigner integrals}}}\\~\\
Solesne Bourguin\footnote{
Faculté des Sciences, de la Technologie et de la Communication; UR en Mathématiques. 6, rue Richard Coudenhove-Kalergi, L-1359 Luxembourg.
Email: {\tt solesne.bourguin@uni.lu}
}
and 
Jean-Christophe Breton\footnote{
IRMAR, CNRS 6625, Université de Rennes 1, 263 avenue du Général Leclerc, 35042 Rennes cedex, France. 
Email: {\tt jean-christophe.breton@univ-rennes1.fr}
}

\end{center}
{\small \noindent {\bf Abstract:} 
We investigate generalizations of the Cramér theorem. 
This theorem asserts that a Gaussian random variable can be decomposed into the sum of independent random variables if and only if they are Gaussian.  
We prove asymptotic counterparts of such decomposition results for multiple Wiener integrals 
and prove that similar results are true for the (asymptotic) decomposition of the semicircular distribution into free multiple Wigner integrals.
\\


\noindent
{\bf 2010 AMS Classification Numbers:} 
60F05, 
60H05, 
46L54. 
\\


\section{Introduction}
\label{sec:intro}
\noindent The Cram\'er theorem states that if $X$ and $Y$ are independent random variables then 
$X+Y$ is Gaussian if and only if $X$ and $Y$ are Gaussian. 
While the reciprocal sense is immediate, the direct sense is due to Cram\'er in \cite{Cramer36} in 1936 (from a conjecture by Lévy). 
Just after, it was shown in 1937 by Raikov that, roughly speaking, the same is true if we consider the class of Poissonian distributions instead of the class of Gaussian distributions.  
We refer to \cite{Lukacs} for a general discussion about the Cram\'er theorem and decomposition of distributions. 

\medskip\noindent
Recently, (asymptotic) counterparts of this result have been investigated by Tudor in \cite{Tudor11} on the Wiener space using Malliavin calculus, see \eqref{eq:Cramer2} for details. 
Using similar arguments in \cite{BT11}, Bourguin and Tudor have explored non-asymptotic and asymptotic Cramér type result for the decomposition of Gamma random variables in terms of independent multiple Wiener integrals.
In this case, the (asymptotic) decompositions obtained in \cite{BT11} are only valid for Wiener integrals, ie. for $X_n=I_{q_1}^W(f_n)$ and $Y_n=I_{q_2}^W(g_n)$.

\medskip\noindent
In this note, we explore similar asymptotic Cram\'er type decompositions for stochastic integrals. 
First, we are interested in Wiener multiple integrals and we obtain an asymptotic decomposition of a Gaussian distribution by such integrals, 
see Proposition \ref{prop:Cramer_Wiener2}. 
This recovers a special case of \cite{Tudor11} but with a short new proof based on the behavior of the moments of multiple Wiener integrals.   
Next, we investigate analogues of the previous result in a free probability setting. 
A motivation for this study is that it is known that several classical probabilistic results have a counterpart in a free probability context. 
A natural question is thus whether there exists a free analogue of the Cram\'er theorem for the semicircular distributions which are the free analogues of Gaussian distributions. 
The question is negatively answered by Bercovici and Voiculescu in \cite{BV95}. 
It is also shown that the Raikov theorem does not have a counterpart in the free probability setting neither, see \cite{BG05}. 
The situation is thus much more complicated for the Cramér theorem in the free probability setting. 
Since in the classical setting we have a short simple proof for multiple Wiener integrals, we investigate the analogous situation in a free setting 
to obtain a positive Cramér type result in a free probability space. 
Recall that, based on the (standard) semicircular distribution, a free Brownian motion can be defined 
and the so-called Wigner integrals are next constructed as the (multiple) stochastic integrals with respect to this free Brownian motion.  
Such integrals are thus the free counterparts of the multiple Wiener integrals. 
Moreover, these Wigner integrals were recently analyzed in \cite{KNPS10} and \cite{NR12} where the required properties were derived.  
In the sequel, we investigate more precisely the asymptotic decomposition of the semicircular distributions by multiple Wigner integrals and we propose Proposition \ref{prop:Cramer_Wigner2} as a first positive result in this direction.  

\medskip\noindent
The paper is organized as follows: 
in Section \ref{sec:tools}, we give the main notations and properties required in the sequel about Wiener integrals, free probability and Wigner integrals.  
We discuss the Cramér theorem and its extensions in Section \ref{sec:Cramer}. 
This section also contains the contribution of the paper for the multiple Wiener integrals (Section \ref{sec:Wiener}) and the multiple Wigner integrals (Section~\ref{sec:Wigner}). 


\section{Probabilistic tools}
\label{sec:tools}
For the sake of self-containedness, we shortly describe in this section the main objects we deal with. 


\subsection{Multiple Wiener integrals}
Let $(W_{t})_{t\geq 0}$ be a classical Wiener process on a standard Wiener space $\left(\Omega,{\mathcal{F}},\pp\right)$. 
We recall that if $f\in L^{2}(\rit_+^{q})$ with $q\geq 1$ integer, $I_q^W(f)$ stands for the multiple Wiener integral of $f$ with respect to $W$. The set $\HH_{q} $ of such integrals is called the $q^{\mbox{\tiny{th}}}$ Wiener chaos. 
The multiple integrals are centered and satisfie the Itô isometry: for  $f\in L^{2}(\rit_+^{q_1})$ and $g\in L^{2}(\rit_+^{q_2})$, $\ee\left[I_{q_1}^W(f) I_{q_2}^W(g)\right]=\ind_{\{q_1=q_2\}} \ q_1! \langle f,g\rangle_{L^2(\rit_+^{q_1})}$.
We refer to basic references such as \cite{Nu06} for details. 
\\
Since the Cramér theorem deals with independent random variables, recall that the independence of two 
Wiener integrals is characterized by its contraction: 
namely, for $q_1,q_2 \geq 1$ and $f\in L^{2}(\rit_+^{q_1})$, $g\in L^{2}(\rit_+^{q_2})$ symmetric functions, $I_{q_1}^{W}(f)$ and $I_{q_2}^{W}(g)$ are independent if and only if $\|f \otimes_1 g \|_{L^2\big(\rit_+^{q_1+q_2-2}\big)} =0$, see \cite{UZ}. 
Above, $f\otimes_1 g$ is an instance of the contraction between $f$ and $g$ which, in general, is defined for $0 \leq \ell \leq q_1 \wedge q_2:=\min(q_1,q_2)$, by 
\begin{eqnarray*}
\nonumber
&&  (f\otimes_{\ell} g) ( t_{1}, \ldots, t_{q_1 + q_2-2\ell})\\
\label{contra}
&&  =\int_{\rit_+^{\ell} } f( t_{1}, \ldots, t_{q_1 - \ell}, s_{1}, \ldots, s_{\ell})
g(t_{q_1 - \ell + 1}, \ldots, t_{q_1 + q_2-2\ell},s_{1}, \ldots,s_{\ell}) \ ds_{1}\ldots ds_{\ell}. 
\end{eqnarray*}
The following theorem gathers several results on the convergence of multiple Wiener integrals. 
It provides powerful criteria in order to prove that a sequence of multiple Wiener integrals converges to the standard normal law 
(see in particular \cite{NP09, NOL08, NP05}).
\begin{theo}
\label{fourthmomentthm}
Fix $q\geq 2$, and consider a sequence $\{F_n : n\geq 1\}$ such
that $F_n=I_q(f_n)$, $n\geq 1$, where $f_n\in L_{s}^{2}(\rit_+^{q})$.
Assume moreover that $\ee[F_n^2]=q!\|f_n\|^2_{L^{2}(\rit_+^{q})}\rightarrow 1$. 
Then, the following two conditions are equivalent, as  $n\to+\infty$:
\begin{enumerate}
\item[\rm(i)] $F_n$ converges in distribution to $Z\sim \mathcal{N}(0,1)$;
\item[\rm(ii)] $\ee[F_n^4]\rightarrow 3$;
\end{enumerate}
\end{theo}
\begin{Rem}
\label{rem:M4}
{\rm 
In complement to (ii) above, observe that for such random variables $F_n$, we have $\ee[F_{n}^4]>3$, see Remark 1.3, point 5 in \cite{BBNP}.
}
\end{Rem}


\subsection{Free probability and Wigner integrals}

In this section, we recall the basic notions and objects of free probability, also called non--commutative probability: 
we revisite, for the unfamiliar reader with this context, the notions of random variable, probability distribution, convergence in law, freeness, semicircular distribution, free Brownian motion and Wigner multiple integral. 
The reader familiar with these notion may skip this section 
while the others may consult \cite{HP00} or \cite{NS06} for a deeper insight.


\medskip\noindent
{\bf Non-commutative probability spaces.}
A free tracial (non--commutative) probability space is a pair $(\mathscr{A},\varphi)$, where $\mathscr{A}$ is a von Neumann algebra (ie. a $\star$--algebra of bounded operators on a  
separable Hilbert space that is closed in the weak operator topology and contains the identity operator) and $\varphi$ is a trace operator (ie. a unital linear functional which is 
weakly continuous, 
positive ($\varphi(X) \geq 0$ whenever $X \in \mathscr{A}$ is non-negative), 
faithful ($\varphi(YY^\ast) = 0 \Rightarrow Y = 0$) 
and tracial ($\varphi(XY) = \varphi(YX)$ for all $X,Y \in \mathscr{A}$)). In this context, we will refer to the self-adjoint elements of the von Neumann algebra $\mathscr{A}$ as random variables.
\\
A standard probability space $(\Omega,{\cal F},\pp)$ is a special case of such non--commutative probability space $(\mathscr{A},\varphi)$: take $\mathscr{A}=L^\infty(\Omega, {\cal F},\pp)$ with the involution $X\mapsto X^*$ given by the complex conjugation and $\varphi(X)=\ee[X]$. 
Another simple example, truly non-commutative, is the space of $(d\times d)$-matrix valued random variables $\mathscr{A}=L^\infty(\Omega, {\cal F},\pp; M_d(\cit))$
where $X^*$ is the usual adjoint of the matrix $X$ and $\varphi(X)=\frac 1d\ee[\mbox{Tr}(X)]$.
\\
Let $(\mathscr{A},\varphi)$ be a free tracial probability space. The law of a random variable $X$ on this space is defined as the unique measure $\mu_{X}$ on the real line such that 
$\int_{\mathbb{R}}P(t)\ d\mu_{X}(t) = \varphi(P(X))$
for all polynomials $P \in \mathbb{R}\left[ X\right]$. We shall note $X\loifree \mu_X$. 
Furthermore, the $k^{\mbox{\tiny{th}}}$ moment of $X$ is the real number $m_{k}(X) = \varphi(X^{k})$. Observe that, by linearity, the probability law of $X$ is determined by its moments.
A sequence $(X_{n})_{n \geq 1}$ of non--commutative random variables is said to converge in law to a limiting random variable $X_{\infty}$ if and only if we have the convergence in the sense of moments, that is 
$\lim_{n \to+\infty}\varphi(P(X_{n})) = \varphi(P(X_{\infty}))$ for any real valued polynomial $P$.
In this case, we write $X_{n} \cvfree X_{\infty}$.
\\~\\
Another central concept in non--commutative probability is \textit{freeness}. It is the counterpart of independence in a classical probability setting. Consider $(\mathscr{A},\varphi)$ a free tracial probability space, and let $\mathscr{A}_{1},\ldots,\mathscr{A}_{p}$ be unital subalgebras of $\mathscr{A}$. 
Let $X_{1},\ldots,X_{m}$ be such that for each $1\leq j<m$, $X_{j} \in \mathscr{A}_{i(j)}$, where $i(1) \neq i(2)$, $i(2) \neq i(3)$, $\ldots$, $i(m-1) \neq i(m)$ and such that $\varphi(X_{j}) = 0$ for all $j$. 
Then, the family of subalgebras $\left\lbrace \mathscr{A}_{1},\ldots,\mathscr{A}_{p} \right\rbrace $ is called free or freely independent if $\varphi(X_{1}X_{2}\cdots X_{m}) = 0$. Random variables are termed free or freely independent if the unital algebras they generate are free. Freeness is a much more complicated concept than classical independence; for example, if $X$ and $Y$ are free, we have $\varphi(XYXY) = \varphi(Y)^{2}\varphi(X^{2}) + \varphi(X)^{2}\varphi(Y^{2}) - \varphi(X)^{2}\varphi(Y)^{2}$, which contrasts to $\ee[XYXY]=\ee[X^2]\ee[Y^2]$ in a classical probability setting.


\medskip\noindent
{\bf Semicircular distributions, free Brownian motion and Wigner integrals.}
The semicircular distribution $\mathcal{S}(m,\sigma^{2})$ with mean $m$ and variance $\sigma^{2} > 0$ is the probability distribution given by the density
\begin{equation*}
\mathcal{S}(m,\sigma^{2})(dx) = \frac{1}{2\pi\sigma^{2}}\sqrt{4\sigma^{2} - (x-m)^{2}}\ \ind_{\left\lbrace \vert x-m \vert \leq 2\sigma\right\rbrace }\ dx.
\end{equation*}
The semicircular distributions play the same role as the Gaussian distributions in the classical setting. 
They enjoy similar properties, for instance if $X \loifree  \mathcal{S}(m_{X},\sigma_{X}^{2})$ and $a,b \in \mathbb{R}$ then $aX + b \loifree  \mathcal{S}(am_{X} + b,a^{2}\sigma^{2})$ ; 
moreover if $Y \loifree  \mathcal{S}(m_{Y},\sigma_{Y}^{2})$ is freely independent from $X$ then 
$X+Y \loifree  \mathcal{S}(m_{X} + m_{Y},\sigma_{X}^{2} + \sigma_{Y}^{2})$.
Note that for $X\loifree{\cal S}(0,1)$, we have $m_1(X)=0$, $m_2(X)=1$ and $m_4(X)=2$. 
\\~\\
A free Brownian motion $(S(t))_{t\geq 0}$ is a non--commutative stochastic process. 
It is a family of self-adjoint operators on a tracial probability space with the following characteristic properties:
\begin{enumerate}
\item $S(0) = 0$.
\item For $0\leq t_1\leq t_2$, the law of $S(t_{2}) - S(t_{1})$ is the semicircular distribution ${\cal S}(0, t_2-t_1)$.
\item For all $n$ and $0<t_{1}< t_{2}<\cdots <t_{n}$, 
the increments $S(t_{1})$, $S(t_{2}) - S(t_{1})$, $\ldots$, $S(t_{n}) - S(t_{n-1})$ are freely independent.
\end{enumerate}
Multiple Wigner integrals of order $q\geq 1$, $I_q^S(f)$, are defined for complex-valued functions $f$ of $L^{2}(\mathbb{R}_{+}^{q})$. They are centered $\varphi\left(I_{q_1}^S(f)\right) = 0$, and the counterpart of the Itô isometry holds true: for $f\in L^{2}(\mathbb{R}_{+}^{q_1}), g\in L^{2}(\mathbb{R}_{+}^{q_2})$, we have 
$\varphi \left(I_{q_1}^S(f)I_{q_2}^S(g) \right) = \ind_{\left\lbrace q_1 = q_2\right\rbrace }\langle f,g^*\rangle_{L^2(\mathbb{R}_{+}^{q_1})}$
where $g^*(t_1,\ldots,t_q)=\overline{g(t_q,\ldots,t_1)}$ and the bar stands for the complex conjugation (see for instance (1.6) in \cite{KNPS10}). The operator $I_q^S$ is an isometry from $L^2(\mathbb{R}_{+}^q)$ to the Hilbert space of operators generated by the free Brownian motion equipped with the inner product $\langle X,Y\rangle_\varphi=\varphi(Y^*X)$ to $L^2(\Omega)$. In order for $I_q^S(f)$ to be a random variable in the sense of the above Definition 
(ie. a self-adjoint operator), it is necessary for $f$ to be {\it mirror--symmetric}, that is, $f=f^*$. 
\\
We also define the nested contractions that will play in the free setting a similar role to the (regular) contractions in the classical setting. The nested contraction $(f\stackrel{\ell}{\frown} g)$ is defined for $0 \leq \ell \leq q_1 \wedge q_2:=\min(q_1,q_2)$ by
\begin{eqnarray*}
\nonumber
&&  (f\stackrel{\ell}{\frown} g) ( t_{1}, \ldots, t_{q_1 + q_2-2\ell})\\
\label{contra2}
&&  =\int_{\rit_+^{\ell} } f( t_{1}, \ldots, t_{q_1 - \ell}, s_{1}, \ldots, s_{\ell})
g(s_{\ell}, \ldots, s_{1}, t_{q_1 - \ell + 1}, \ldots, t_{q_1 + q_2-2\ell}) \ ds_{1}\ldots ds_{\ell}. 
\end{eqnarray*}
The following result by Nourdin and Rosinski, see \cite{NR12}, shows that like for the multiple Wiener integrals,  when the algebra ${\cal A}$ is real, the (nested) contractions characterize the freeness of multiple Wigner integrals. 
\begin{prop}
\label{prop:Nourdin_indep}
Consider a real von Neumann algebra ${\cal A}$. 
Let $q_1,q_2 \geq 1$, and let real-valued $f\in L^{2}(\rit_+^{q_1})$ and $g\in L^{2}(\rit_+^{q_2})$ be fully symmetric functions. 
Then, $I_{q_1}^{S}(f)$ and $I_{q_2}^{S}(g)$ are freely independent if and only if 
$\Vert f \stackrel{1}{\frown} g \Vert _{L^{2}(\rit_+^{q_1+q_2-2} )} =0$.
Actually, $I_{q_1}^{S}(f)$ and $I_{q_2}^{S}(g)$ are freely independent if and only if $I_{q_1}^{W}(f)$ and $I_{q_2}^{W}(g)$ are classically independent.
\end{prop}

\noindent
The similarities between the Wiener and the Wigner integrals lead Kemp, Nourdin, Peccati and Speicher 
to prove a \textit{transfer principle} between the classical and the free chaoses (see \cite{KNPS10}).
\begin{prop}[Transfer principle] 
\label{prop:transfer}
Let $q\geq 2$ be an integer, and let $(f_n)_{n\in\nit}$ be a sequence of fully symmetric functions in $L^2(\rit_+^q)$. 
Let $\sigma>0$ be a finite constant. 
Then, as $n\to+\infty$. 
\begin{enumerate}
\item $\lim_{n\to+\infty}\ee[I_q^W(f_n)^2]=q!\sigma^2$ if and only if $\lim_{n \rightarrow+\infty}\varphi(I_q^S(f_n)^2)=\sigma^2$.
\item If the asymptotic relations in 1. are verified, then $I_q^W(f_n) \cvlaw \mathcal{N}(0,q!\sigma^2)$ 
if and only if $I_q^S(f_n) \cvfree S(0,\sigma^2)$.
\end{enumerate}
\end{prop}
In the sequel, the first moments of the multiple integrals have some importance.
For this purpose, we summarize the following useful results from \cite{KNPS10}: 
1. is given by Corollary 1.7 while 2. is the main result, Theorem 1.3, therein. 
\begin{prop}[Moments of a multiple Wigner integrals]
\label{prop:Wigner_moment}
$ $
\begin{enumerate}
\item For $f\in L^2(\rit_+^q)$ not identically equal to zero, we have $\varphi(I_q^{S}(f)^4)>2(\varphi(I_q^S(f)^2)^2$.
\item Let $q\geq 2$ be an integer and let $f_n\in L^2(\rit_+^q)$ be a sequence of mirror-symmetric functions with $\|f_n\|_{L^2(\rit_+^q)}=1$. 
Then, when $n\to+\infty$,  
$\varphi(I_q^S(f_n)^4)\to 2$ {\it iff} $I_q^S(f_n)\cvfree {\cal S}(0,1)$.
\end{enumerate}
\end{prop}
As a consequence of 1., it comes that a non--zero multiple Wigner integral cannot have a semicircular distribution 
(recall that $m_4(X)=2$ for $X\loifree {\cal S}(0,1)$). 


\section{Cramér type theorems}
\label{sec:Cramer}


\subsection{Classical results}
The seminal result motivating this line of work is the following:
\begin{theo}[Cramér, 1936]
\noindent 
Let $X$ and $Y$ be two centered independent real valued random variables and let $Z = X+Y$. 
Then, it holds true that
\begin{equation}
\label{eq:Cramer}
\Big\{Z \sim \mathcal{N}(0,\sigma_{1}^{2} + \sigma_{2}^{2})\Big\} 
\Longleftrightarrow 
\Big\{X \sim \mathcal{N}(0,\sigma_{1}^{2})  \mbox{\ \ and\ \ }  Y \sim \mathcal{N}(0,\sigma_{2}^{2})\Big\}.
\end{equation}
\end{theo}
One direction is elementary since the sum of two independent Gaussian random variables is necessarily Gaussian. 
The second direction is less trivial. 
It was conjectured by Lévy and proved by Cramér in \cite{Cramer36} with powerful results from complex analysis. 
A similar result was proved just after by Raikov in \cite{Raikov37} (see also \cite{Raikov38}) for the class of Poissonian distributions $\{{\cal P}(\lambda), \lambda>0\}$. 
\\~\\
A natural generalization of Cram\'er's result consists in obtaining asymptotic counterparts of (\ref{eq:Cramer}) where the equalities in law are replaced by convergences in law, namely with obvious notations and assuming $X_n$ and $Y_n$ are independent for all $n\geq 1$
\begin{equation}
\label{eq:Cramer2}
\Big\{Z_n \Rightarrow \mathcal{N}(0,\sigma_{1}^{2} + \sigma_{2}^{2})\Big\} 
\quad \Longleftrightarrow \quad 
\Big\{X_n \Rightarrow \mathcal{N}(0,\sigma_{1}^{2})  \mbox{\ \ and\ \ }  Y_n\Rightarrow \mathcal{N}(0,\sigma_{2}^{2})\Big\},
\end{equation}
where, here and hereafter, the limit are taken with $n\to+\infty$. 
Such a generalization has recently been investigated by Tudor in \cite{Tudor11} for square integrable random variables on the Wiener space using Malliavin calculus. More recently, with similar arguments, Bourguin and Tudor have shown in \cite{BT11} that the asymptotic Cramér type decomposition holds true for multiple Wiener integrals converging to Gamma distributions.
\\~\\
In Section \ref{sec:Wiener}, we recover a special case of the asymptotic result by Tudor (see \cite{Tudor11}) for multiple Wiener integrals with a short new proof of the asymptotic Cram\'er equivalence. The proof is independent from Cramér's original result and is based on the behavior of the first moments of the multiple Wiener integrals. 
\\~\\
Another generalization of Cram\'er's result would be to derive such an analogous in a free probability context by replacing the concept of independence by freeness, but the situation is much more complicated in the free probability context 
and it is shown by Bercovici and Voiculescu that no such general result holds true (see \cite{BV95} for details). The difficulties in the free probability context are confirmed 
in \cite{BG05} where the failure of (a counterpart to) the 
Raikov theorem is shown for free Poisson distributions (the so-called Marchenko-Pastur).
\\~\\
This is thus a non trivial problem to obtain positive result for Cramér type decompositions in a free probability setting.  
Since similarities have been observed between Wiener and Wigner integrals, as illustrated for instance by the transfer principle, 
and since Cramér's decompositions are known for multiple Wiener integrals, 
we address (in Section \ref{sec:Wigner}) the problem of Cramér type decompositions for free random variables having the form of mutiple Wigner integrals. This is the object of Proposition \ref{prop:Cramer_Wigner2}.


\subsection{Wiener integrals}
\label{sec:Wiener}

In this section, we investigate Cramér type decomposition for multiple Wiener integrals. 
First, since the Wiener chaoses of order $q\geq 2$ do not contain Gaussian distributions, 
the non-asymptotic result (\ref{eq:Cramer}) is irrelevant for multiple Wiener integrals of order $q\geq 2$ with positive variance.   
We are thus interested in asymptotic decomposition (\ref{eq:Cramer}) and we prove (for multiple Wiener integrals) a stronger result than the one contained in \cite{Tudor11}. 
The proof we propose below in this context is new, short and independent from Cramér's original result: it is based on a recent result by Nourdin and Poly on convergence in total variation on Wiener chaoses (see \cite{NP12}) as well as on the fact that on the Wiener chaoses, the central convergence is controlled by the convergence of the moments of order $2$ and $4$ (see \cite{NOL08, NP05}). 
Such a phenomenon has been recently furtherly investigated in \cite{BBNP} where optimal Berry-Ess\'een rates are given in terms of third and fourth cumulants. 
Moreover, the proof is easily adapted in Section \ref{sec:Wigner} to derive a similar result for Wigner integrals and semicircular distribution. 

\begin{prop}
\label{prop:Cramer_Wiener2}
Fix $q_{1}, q_{2} \geq 1$ and let $(f_n)_{n\in\nit}$ (resp. $(g_n)_{n\in\nit}$) be a sequence of symmetric functions in $L^{2}(\mathbb{R}_{+}^{q_{1}})$ (resp. in $L^{2}(\mathbb{R}_{+}^{q_{2}})$). Set $X_{n} = I_{q_1}^{W}(f_n)$ and $Y_n = I_{q_2}^{W}(g_n)$. Assume that $\ee\left[X_n^3 Y_n \right] = \ee\left[X_n Y_n^3 \right] = \ee\left[X_n Y_n \right] = 0$ for any $n$ and that the sum $X_n + Y_n \cvlaw \mathcal{N}\left(0,\sigma^2 \right)$, $\sigma^2 > 0$. Then the sequence $(X_n,Y_n)$ is tight. Moreover, if $(n_k)_{k\geq 1}$ is a subsequence such that $(X_{n_k},Y_{n_k}) \cvlaw (X_{\infty},Y_{\infty})$ then, necessarily, $(X_{\infty},Y_{\infty}) \sim \mathcal{N}_{2}\left(0,\rm{diag}\left( \sigma_{1}^{2},\sigma_{2}^{2}\right)  \right) $ with $\sigma_{1}^{2} + \sigma_{2}^{2} = \sigma^{2}$.  
\end{prop}

\medskip\noindent
\begin{demo}
Only the direct implication requires an argument. Since $X_n + Y_n$ is tight, Lemma 2.4 in \cite{NP12} implies that $X_n + Y_n$ is uniformly bounded in all the $L^p$. But as $\ee\left[X_n Y_n \right] = 0$, both $X_n$ and $Y_n$ are uniformly bounded in $L^2$, and in all the $L^p$ as well by hypercontractivity. 
This implies that the sequence $(X_n,Y_n)$ is tight. 

\medskip\noindent
Now, let $(n_k)_{k\geq 1}$ be a subsequence such that $(X_{n_k},Y_{n_k}) \cvlaw (X_{\infty},Y_{\infty})$. Set $Z_{n_k}=X_{n_k}+Y_{n_k}$, $\lim_{n \rightarrow +\infty}\Var(X_n)=\lim_{n\rightarrow +\infty}q_1!\|f_n\|_{L^2(\rit_+^{q_1})}^2 = \sigma_1^2$, $\lim_{n \rightarrow +\infty}\Var(Y_n)=\lim_{n\to+\infty}q_2!\|g_n\|_{L^2(\rit_+^{q_2})}^2 = \sigma_2^2$ and $N \sim \mathcal{N}(0,1)$. Without loss of generality, assume that $\sigma_1^2+\sigma_2^2=1$.
Because $X_n$ and $Y_n$ are centered sequences and because $\ee\left[X_n Y_n \right] = 0$, we have 
\begin{equation}
\label{eq:11}
q_1!\|f_{n_k}\|_{L^2(\rit_+^{q_1})}^2+q_2!\|g_{n_k}\|_{L^2(\rit_+^{q_2})}^2
=\Var(X_{n_k})+\Var(Y_{n_k})=\Var(Z_{n_k}) \underset{k \rightarrow +\infty}{\longrightarrow} 1,
\end{equation}
and necessarily $\sigma^2=1$. 
The hypercontractivity property of the Wiener chaoses implies that $\ee[Z_{n_k}^p] \underset{k \rightarrow +\infty}{\longrightarrow} \ee[N^p]$ for every integer $p\geq 3$
(see \cite[Chap. 5]{Janson} or Remark 1.8 in \cite{BBNP}). 
In particular, one has $\ee[Z_{n_k}^3] \underset{k \rightarrow +\infty}{\longrightarrow} 0$ and $\ee[Z_{n_k}^4] \underset{k \rightarrow +\infty}{\longrightarrow} 3$. Combining this with the fact that $\ee\left[X_n^3 Y_n \right] = \ee\left[X_n Y_n^3 \right] = 0$ for any $n$ and using \eqref{eq:11}, we can write
\begin{eqnarray}
\nonumber
&&\lim_{n\to +\infty}\ee[Z_{n_k}^4]-3\\
\nonumber
&=&\lim_{n\to +\infty}\left(\ee[X_{n_k}^4]+\ee[Y_{n_k}^4]+6\ee[X_{n_k}^2 Y_{n_k}^2]-3\left(q_1!\|f_{n_k}\|_{L^2(\rit_+^{q_1})}^2+q_2!\|g_{n_k}\|_{L^2(\rit_+^{q_2})}^2\right)^2\right) \\
\nonumber
&=&\lim_{n\to +\infty}\Big(\left(\ee[X_{n_k}^4]-3(q_1!)^2\|f_{n_k}\|_{L^2(\rit_+^{q_1})}^4\right)
+\left(\ee[Y_{n_k}^4]-3(q_2!)^2\|g_{n_k}\|_{L^2(\rit_+^{q_2})}^4\right)\\
\label{eq:00}
&&\hskip 1.4cm +6\left(\ee[X_{n_k}^2 Y_{n_k}^2]-q_1!q_2!\|f_{n_k}\|_{L^2(\rit_+^{q_1})}^2\|g_{n_k}\|_{L^2(\rit_+^{q_2})}^2\right)\Big).  
\end{eqnarray} 
Now recall that for $F$ belonging to a fixed Wiener chaos, with unit variance, one has $\ee[F^4]>3$ (see Remark \ref{rem:M4}). This implies that $\ee\left[ \left( \frac{X_{n_k}}{\sqrt{q_1!}\|f_{n_k}\|_{L^2(\rit_+^{q_1})}}\right) ^4\right] > 3$ and that $\ee[X_{n_k}^4]-3(q_1!)^2\|f_{n_k}\|_{L^2(\rit_+^{q_1})}^4 > 0$. Similarly, $\ee[Y_{n_k}^4]-3(q_2!)^2\|g_{n_k}\|_{L^2(\rit_+^{q_2})}^4 > 0$. Because $\underset{k \rightarrow \infty}{\rm lim}\ee[Z_{n_k}^4]-3 = 0$, it follows that the term in (\ref{eq:00}) converges to zero. As a consequence, 
$$
\ee[X_{n_k}^4]-3(q_1!)^2\|f_{n_k}\|_{L^2(\rit_+^{q_1})}^4 \underset{k \rightarrow +\infty}{\rightarrow} 0
\quad \mbox{ and } \quad 
\ee[Y_{n_k}^4]-3(q_2!)^2\|g_{n_k}\|_{L^2(\rit_+^{q_2})}^4 \underset{k \rightarrow +\infty}{\rightarrow} 0. 
$$
Applying Theorem \ref{fourthmomentthm} to $\widetilde X_{n_k}=\frac{X_{n_k}}{\sqrt{q_1!} \|f_{n_k}\|_{L^2(\rit_+^{q_1})}}$ 
and 
$\widetilde Y_{n_k}=\frac{Y_{n_k}}{\sqrt{q_2!} \|g_{n_k}\|_{L^2(\rit_+^{q_2})}}$ yields  
$$
X_{n_k} \cvlaw \mathcal{N}(0,\sigma_{1}^{2})
\quad \mbox{ and } \quad 
Y_{n_k} \cvlaw \mathcal{N}(0,\sigma_{2}^{2}). 
$$
Because $\ee\left[X_n Y_n \right] = 0$, Peccati and Tudor's convergence theorem for vector-valued multiple stochastic integrals (see \cite{PT05}, Theorem 1) concludes the proof.
\CQFD
\end{demo}
\begin{Rem}
{\rm 
\begin{itemize}
\item The moment conditions of Proposition \ref{prop:Cramer_Wiener2} are satisfied in particular if $X_n$ and $Y_n$ are independent for all $n$. But this is obviously a much stronger requirement. 
\item If $\sigma_1>0$ and $\sigma_2>0$ then the convergence to $\mathcal{N}_{2}\left(0,\rm{diag}\left( \sigma_{1}^{2},\sigma_{2}^{2}\right)  \right) $
holds in total variation, see Theorem 5.2 in \cite{NP12}. 
\item Proposition \ref{prop:Cramer_Wiener2} can easily be generalized to sums of independent Wiener integrals without any further argument. 
\end{itemize}
}
\end{Rem}


\subsection{Wigner integrals}
\label{sec:Wigner}

In this section, we investigate a free analogue of the Cramér type decomposition. We deal with multiple Wigner integrals and give similar free results as in Section~\ref{sec:Wiener}. 
Recall that the Cramér theorem is not true in general for the semicircular distribution, see \cite{BV95}, 
and it is in particular not true for decomposition into multiple Wigner integrals $I_q^{S}(f)$ since such an integral cannot be semicircular (see the consequence of Proposition \ref{prop:Wigner_moment}).
Our goal is thus to derive an asymptotic Cramér type decomposition result for multiple Wigner integrals. 

\medskip\noindent
First, in the real case, for freely independent multiple Wigner integrals $X_n=I_q^S(f_n)$ and $Y_n=I_q^S(g_n)$ of the same order $q$, 
when $f_n, g_n$ are real-valued fully symmetric kernels, the desired result is easily obtained from Proposition~\ref{prop:Cramer_Wiener2} by the transfer principle (Proposition \ref{prop:transfer}). 
In order to illustrate the use of this transfer principle, we give here a proof of this fact where, for that purpose, we assume that the von Neumann algebra ${\cal A}$ is real. 

\medskip\noindent
Let $Z_n=X_n+Y_n=I_q^S(f_n+g_n) \cvfree S(0,\sigma^2)$ with $X_n$ and $Y_n$ freely independent 
with $\|f_n\|_{L^2(\rit_+^q)}^2\to \sigma_1^2$ and $\|g_n\|_{L^2(\rit_+^q)}^2\to \sigma_2^2$. 
Recall that all the limits are taken with $n\to+\infty$.
\\  
Since $I_q^S(f_n)$ and $I_q^S(g_n)$ are freely independent and $f_n$ and $g_n$ are fully symmetric, 
we have $f_n\stackrel{1}{\frown}g_n=0$ by Proposition \ref{prop:Nourdin_indep}. 
As a consequence $\langle f_n,g_n\rangle_{L^2(\rit_+^q)}=f_n\stackrel{q}{\frown}g_n=0$ 
and
\begin{eqnarray*}
\|f_n+g_n\|_{L^2(\rit_+^q)}^2&=&
\|f_n\|_{L^2(\rit_+^q)}^2+\|g_n\|_{L^2(\rit_+^q)}^2+2\langle f_n,g_n\rangle_{L^2(\rit_+^q)}\\
&=&\|f_n\|_{L^2(\rit_+^q)}^2+\|g_n\|_{L^2(\rit_+^q)}^2\\
&\to& \sigma_1^2+\sigma_2^2. 
\end{eqnarray*}
Since $\varphi(I_q^S(f_n+g_n)^2)\to \sigma_1^2+\sigma_2^2$, Proposition \ref{prop:transfer} entails 
$$
I_q^W(f_n+g_n)=I_q^W(f_n)+I_q^W(g_n)
\Longrightarrow {\mathcal N}\big(0,q!(\sigma_1^2+\sigma_2^2)\big). 
$$ 
Since by the Üstünel-Zakai criterion  $I_q^W(f_n)$ and $I_q^W(g_n)$ are (classically) independent, see \cite{UZ}, 
Proposition \ref{prop:Cramer_Wiener2} applies and gives that 
$I_q^W(f_n)\Rightarrow {\mathcal N}(0,q!\sigma_1^2)$ and $I_q^W(g_n)\Rightarrow {\mathcal N}(0,q!\sigma_2^2)$. 
Other two applications of the transfer principle (Proposition \ref{prop:transfer}) give
$I_q^S(f_n)\cvfree S(0,\sigma_1^2)$ and $I_q^S(g_n)\cvfree S(0,\sigma_2^2)$ 
and necessarily $\sigma^2=\sigma_1^2+\sigma_2^2$. 
\CQFD
 
\medskip\noindent
Actually, we obtain a better result for general (complex) von Neumann algebra ${\cal A}$, valid for Wigner integrals of arbitrary (different) orders and with (only) mirror-symmetric kernels, 
by following the same proof as in Proposition~\ref{prop:Cramer_Wiener2}.  
In particular, the key result is the main result in \cite{KNPS10}, Theorem~1.3. 
\begin{prop}
\label{prop:Cramer_Wigner2}
Let $X_n=I_{q_1}^S(f_n)$ and $Y_n=I_{q_2}^S(g_n)$ be free multiple Wigner integrals of orders $q_1$ and $q_2$ with $f_n\in L^2(\rit_+^{q_1})$ and $g_n\in L^2(\rit_+^{q_2})$ mirror symmetric. 
Assume that $\|f_n\|_{L^2(\rit_+^{q_1})}^2\to \sigma_1^2$ and $\|g_n\|_{L^2(\rit_+^{q_2})}^2\to \sigma_2^2$. 
Then, it holds that
\begin{equation}
\label{eq:Cramer4}
\Big\{X_n + Y_n \cvfree \mathcal{S}(0,\sigma^{2})\Big\} 
\quad \Longleftrightarrow \quad 
\Big\{X_n \cvfree \mathcal{S}(0,\sigma_{1}^{2})  \mbox{\ \ and\ \ }  Y_n\cvfree \mathcal{S}(0,\sigma_{2}^{2})\Big\}
\end{equation}
and in this case, necessarily, $\sigma^2=\sigma_1^2+\sigma_2^2$. 
\end{prop}

\noindent
\begin{demo}
The direction $\Leftarrow$ in \eqref{eq:Cramer4} is a straightforward consequence of basic properties of the semicircurlar distributions.
The direction $\Rightarrow$ follows a similar scheme as for Proposition \ref{prop:Cramer_Wiener2}. 
Assume that $X_n+Y_n\cvfree S$ where $S$ stands for a free random variable with distribution $S(0,\sigma^2)$. 
Without loss of generality, we also assume that $\sigma^2=1$. 
By definition of the free weak convergence, we have 
$\varphi(P(X_n+Y_n))\to \varphi(P(S))$ for any polynomial $P$.
A straightforward computation using the tracial property of $\varphi$ as well as the fact that $X_n$ and $Y_n$ are free with $\varphi(X_n)=\varphi(Y_n)=0$ entails
\begin{equation}
\label{eq:freecv1}
\varphi(X_n^2)+\varphi(Y_n^2)=\varphi((X_n+Y_n)^2)\lto\varphi(S^2)=1
\end{equation}
and 
\begin{equation}
\label{eq:freecv2}
\varphi(X_n^4)+\varphi(Y_n^4)+4\varphi(X_n^2)\varphi(Y_n^2)
=\varphi((X_n+Y_n)^4)\to \varphi(S^4)=2. 
\end{equation}
Equations \eqref{eq:freecv1} and \eqref{eq:freecv2} together entail
\begin{eqnarray*}
0&=&\lim_{n \to+\infty}\left(\varphi(X_n^4)+\varphi(Y_n^4)+4\varphi(X_n^2)\varphi(Y_n^2)-2\right)\\
&=&\lim_{n \to+\infty}\left(\varphi(X_n^4)+\varphi(Y_n^4)+4\varphi(X_n^2)\varphi(Y_n^2)-2\left(\varphi(X_n^2)+\varphi(Y_n^2)\right)^2\right)\\
&=&\lim_{n \to+\infty}\left(\varphi(X_n^4)+\varphi(Y_n^4)+4\varphi(X_n^2)\varphi(Y_n^2)-2\varphi(X_n^2)^2-2\varphi(Y_n^2)^2-4\varphi(X_n^2)\varphi(Y_n^2)\right)\\
&=&\lim_{n \to+\infty}\left(\left(\varphi(X_n^4)-2\varphi(X_n^2)^2\right)+\left(\varphi(Y_n^4)-2\varphi(Y_n^2)^2\right)\right). 
\end{eqnarray*}
Since by 1) in Proposition \ref{prop:Wigner_moment}, both summands in the above right-hand side are positive, we have 
$$
\varphi(X_n^4)-2\varphi(X_n^2)^2 \to 0 \quad \mbox{ and } \quad 
\varphi(Y_n^4)-2\varphi(Y_n^2)^2\to 0. 
$$
The criterion for free convergence in law to the semicircular from \cite{KNPS10} stated in 2) of Proposition \ref{prop:Wigner_moment} applies and yields
$$
I_{q_1}^S\left(\frac{f_n}{\|f_n\|_{L^2(\rit_+^{q_1})}}\right) \cvfree S(0,1) \quad\mbox{and} \quad
I_{q_2}^S\left(\frac{g_n}{\|g_n\|_{L^2(\rit_+^{q_2})}}\right) \cvfree S(0,1).
$$
Since $\lim_{n \to+\infty}\|f_n\|_{L^2(\rit_+^{q_1})}=\sigma_1$ and $\lim_{n \to+\infty}\|g_n\|_{L^2(\rit_+^{q_2})}=\sigma_2$, we derive 
$X_n \cvfree S(0,\sigma_1^2)$ and $Y_n \cvfree S(0,\sigma_2^2)$ 
and we obtain necessarily that $\sigma_1^2+\sigma_2^2=1$. 
\CQFD
\end{demo}

\medskip\noindent
We can easily adapt the proof of Proposition \ref{prop:Cramer_Wigner2} to obtain: 
\begin{corol}
\label{corol:Cramer_Wigner2}
Let, for $j\geq 1$, $X_n^j=I_{q_j}^S(f_n^j)$, $n\geq 1$, be a sequence of free multiple Wigner integrals of orders $q_j$ with $f_n^j\in L^2(\rit_+^{q_j})$ mirror symmetric. 
Assume that $X_n^j$, $j\geq 1$, are freely independent 
and that $\|f_n^j\|_{L^2(\rit_+^{q_j})}^2\to \sigma_j^2$. 
Then
$$
\Big\{\sum_{j\geq 1} X_n^j \cvfree \mathcal{S}(0,\sigma^2)\Big\} 
\quad \Longleftrightarrow \quad 
\Big\{X_n^j \cvfree \mathcal{S}(0,\sigma_{j}^{2})  \ \mbox{ for all } j\geq 1\Big\}.
$$
and in this case, necessarily, $\sigma^2=\sum_{j\geq 1}\sigma_j^2$. 
\end{corol}
Since every square integrable random variable in the von Neumann algebra ${\cal A}(S)$ generated by the free Brownian motion $S$ 
has a chaotic expansion into multiple Wigner integrals $F=\sum_{q=0}^{+\infty} I_q^S(f_q)$
(see Proposition 5.3.2 in \cite{BS98}), 
Corollary \ref{corol:Cramer_Wigner2} gives a criterion for the convergence of such random variable with freely independent chaos components  
to the semicircular distribution.

\medskip\noindent
The proof follows the same lines with the following lemma replacing \eqref{eq:freecv1} and \eqref{eq:freecv2}. 
\begin{lemme}
Let $X_1, X_2,\dots, X_p, \dots$ be freely independent random variables with $\varphi(X_i)=0$, $i\geq 1$.
We have 
\begin{eqnarray}
\label{eq:varphi2}
\varphi\left(\big(\sum_{i=1}^{+\infty}X_i\big)^2\right)&=&\sum_{i=1}^{+\infty}\varphi(X_i^2)\\
\label{eq:varphi4}
\varphi\left(\big(\sum_{i=1}^{+\infty} X_i\big)^4\right)&=&\sum_{i=1}^{+\infty} \varphi(X_i^4)+4\sum_{i<j}\varphi(X_i^2)\varphi(X_j^2).
\end{eqnarray}
\end{lemme}
\begin{demo} 
First, for finite sum, the proof proceeds by induction. 
The equality \eqref{eq:varphi2} is initiated by \eqref{eq:freecv1} and if it holds true for $p-1$ then 
\begin{eqnarray*}
\nonumber
\varphi\big((X_1+\dots+X_p)^2\big)&=&\varphi\big((X_1+\dots+ X_{p-1})+X_p)^2\big)
=\varphi\big((X_1+\dots+ X_{p-1})^2\big)+\varphi(X_p^2)\\
\nonumber
&=&
\sum_{i=1}^{p-1}\varphi(X_i^2)+\varphi(X_p^2)
=\sum_{i=1}^p\varphi(X_i^2)
\end{eqnarray*}
using the induction hypothesis. 
This establish \eqref{eq:varphi2} for a sum of $p$ terms. 
Similarly, \eqref{eq:freecv2} is initialized by \eqref{eq:freecv2} and if it holds true for $p-1$, then 
\begin{eqnarray*}
\nonumber
\varphi\big((X_1+\dots +X_p)^4\big)
&=&\varphi\big(((X_1+\dots +X_{p-1})+X_p)^4\big)\\
\nonumber
&=&\varphi\big(((X_1+\dots +X_{p-1})^4\big)+\varphi(X_p)^4+4\varphi\big((X_1+\dots+X_{p-1})^2\big)\varphi(X_p^2)\\
&=&\sum_{i=1}^{p-1} \varphi(X_i^4)+4\sum_{i<j}\varphi(X_i^2)\varphi(X_j^2)+\varphi(X_p)^4+4\sum_{i=1}^p \varphi(X_i^2)\varphi(X_p^2)
\end{eqnarray*}
using the induction hypothesis and \eqref{eq:varphi2}. 
This establish \eqref{eq:varphi4} for a sum of $p$ terms. 
Next, the results extend to infinite sums by continuity of $\varphi$. 
\CQFD
\end{demo}

\noindent As far as we know, it is an open problem whether the Cramér type result \eqref{eq:Cramer4} holds for more general random variables in a free probability context. 
Other technics are required to address this problem. 

\medskip\noindent
{\bf Acknowledgments.}
The authors are grateful to an anonymous referee for his/her thorough review and highly appreciate the comments and 
suggestions, which significantly contributed to improving this paper and especially Proposition \ref{prop:Cramer_Wiener2}.

\end{document}